\documentclass{article}
\usepackage{nips07submit_e,times,psfrag,graphicx}

\title{Convex Algorithms for\\ Nonnegative Matrix Factorization}

\author{
Vijay Krishnamurthy\\
ORFE, Princeton University\\
Princeton, NJ 08544\\
\texttt{kvijay@princeton.edu} \\
\And
Alexandre d'Aspremont\\
ORFE, Princeton University\\
Princeton, NJ 08544\\
\texttt{aspremon@princeton.edu} \\
}

\newcommand{\BEAS}{\begin{eqnarray*}}
\newcommand{\EEAS}{\end{eqnarray*}}
\newcommand{\BEA}{\begin{eqnarray}}
\newcommand{\EEA}{\end{eqnarray}}
\newcommand{\BEQ}{\begin{equation}}
\newcommand{\EEQ}{\end{equation}}
\newcommand{\BIT}{\begin{itemize}}
\newcommand{\EIT}{\end{itemize}}
\newcommand{\BNUM}{\begin{enumerate}}
\newcommand{\ENUM}{\end{enumerate}}

\newcommand{\BA}{\begin{array}}
\newcommand{\EA}{\end{array}}
\newcommand{\BC}{\begin{center}}
\newcommand{\EC}{\end{center}}

\newcommand{\ones}{\mathbf 1}

\newcommand{\reals}{{\mbox{\bf R}}}

\newcommand{\symm}{{\mbox{\bf S}}}  

\newcommand{\Rank}{\mathop{\bf Rank}}
\newcommand{\Tr}{\mathop{\bf Tr}}

\newcommand{\QED}{~~\rule[-1pt]{6pt}{6pt}}

\newtheorem{theorem}{Theorem}

\newtheorem{remark}[theorem]{Remark}

\newcounter{exno}

{\begin{quote}}{\end{quote}}

\makeatletter
\long\def\@makecaption#1#2{
   \vskip 9pt 
   \begin{small}
   \setbox\@tempboxa\hbox{{\bf #1:} #2}
   \ifdim \wd\@tempboxa > 5.5in
        \begin{center}
        \begin{minipage}[t]{5.5in}
        \addtolength{\baselineskip}{-0.95pt}
        {\bf #1:} #2 \par
        \addtolength{\baselineskip}{0.95pt}
        \end{minipage}
        \end{center}
   \else 
    \hbox to\hsize{\hfil\box\@tempboxa\hfil}  
   \fi
   \end{small}\par
}
\makeatother

\newcounter{oursection}

\newcounter{lecture}

\begin{document}

\maketitle

\begin{abstract}
We derive approximation algorithms for the nonnegative matrix factorization problem, i.e. the problem of factorizing a matrix as the product of two matrices with nonnegative coefficients. We form convex approximations of this problem which can be solved efficiently and test our algorithms on some classic numerical examples.
\end{abstract}

\section{Introduction}
Nonnegative matrix factorization (NMF) is a classic unsupervised technique to learn a parts-based representation of the data in an additive setting. As such, it is used as a factor analysis tool for high-dimensional data whose components are constrained to be nonnegative. Given a data matrix $A\in\reals^{m \times n}$, we write the nonnegative matrix factorization problem as:
\BEQ\label{eq:nmf}
\BA{ll}
\mbox{minimize} & \mathbf{loss}(A,UV^T)\\
\mbox{subject to} & U,V \geq 0,\\
\EA\EEQ
in the variables $U\in\reals^{m \times k}$ and $V\in\reals^{n \times k}$, where $\mathbf{loss}(X,Y)$ is a loss function, and $k$ is a given rank target. This apparently simple problem can be traced back to \cite{Paat94} and \cite{Lee99} and has found many applications in machine learning and statistics. It was used for example in gene expression data analysis (see e.g. \cite{Brun04}), in signal processing (see \cite{Sha05}), as a clustering tool (see e.g. \cite{Xu03} and \cite{Lang05a}), for image analysis (see \cite{Dono03}), etc. If $k\geq\min(m,n)$, $A=A\mathbf{I}$ is always a solution, so our objective is to make this representation as parsimonious as possible and keep $k$ small. The decomposition is of course not unique and this and other consistency issues were explored in \cite{Dono03}. As a factor analysis technique, NMF is very similar to Principal Component Analysis, however PCA amounts to a simple singular value decomposition of the data matrix which is computationally easy, while NMF is a NP-Hard problem. Furthermore, in all of these applications, the nonnegativity constraint on the components is the result of some physical property of the data, hence cannot be lifted.

\subsection{Current methods}
Here, we briefly summarize the main types of algorithms currently used to solve it. We refer the reader to \cite{Berr06} for a more complete survey. The algorithms listed below have all been implemented in MATLAB libraries such as NMFLAB by \cite{Zdun06} or NMFPACK by \cite{Hoye04}.
\paragraph{Multiplicative update.}
The original algorithm in \cite{Lee99}, when the loss is given by the Mean Squared Error, updates the current iterates $U$ and $V$ as follows:
\[
V^+_{ij}=V_{ij}\frac{(U^TA)_{ji}}{(U^TUV^T)_{ji}} \quad \mbox{and} \quad 
U^+_{ij}=U_{ij}\frac{(AV)_{ij}}{(UV^TV)_{ij}}.
\]
Similar updates exist for other loss functions. This is a descent method but, in this form, it is not guaranteed to converge to a local minimum.
\paragraph{Gradient descent.}
Another set of algorithms (see \cite{Lee01} or \cite{Hoye04} among others) directly apply a gradient descent algorithm to problem (\ref{eq:nmf}).
\paragraph{Alternating least squares.}
When the loss function is given by the MSE, a third group of algorithms take advantage of the fact that the minimization problem in only one of the variables is equivalent to a least-squares problems. These block-coordinate descent algorithms (see \cite{Heil06} among others) minimize the loss by alternating between the LS problems in $U$ and $V$.

Recently, \cite{Hoye04} (see also \cite{Heil06}) also added a penalty term to problem (\ref{eq:nmf}) to make the matrices $U$ and $V$ sparse, which improves interpretability in imaging applications for example.

\subsection{Contribution}
All the algorithms above have one common characteristic, they all solve a \emph{nonconvex} formulation of the nonnegative matrix factorization problem. In particular, this means that they all seek local solutions to the original problem. This creates \emph{stability} issues, i.e. the solutions are very sensitive to the choice of initial point, it also creates \emph{complexity} issues, meaning that no precise bound can be given on the computing time required to solve the problem and suboptimality cannot be measured by computing the duality gap.
Here, we begin by formulating a convex approximation to nonnegative matrix factorization, we then solve the approximate problem using convex optimization methods. This means that we find \emph{global} (hence potentially more stable) solutions to the approximate problem with guaranteed complexity bounds.

In the symmetric case, we first show that the NMF problem can be formulated as the problem of approximating a given matrix by a completely positive matrix. We then use a convex representation of a restriction of the set of completely positive matrices to write a convex restriction of the symmetric NMF problem. In other words, we show that solving problem (\ref{eq:nmf}) over a subset of all the possible choices of $U$ and $V$ is equivalent to a \emph{convex} problem. We then extend these results to the nonsymmetric case.

The paper is organized as follows, in Section \ref{s:cvx-approx} we detail our convex approximations of problem (\ref{eq:nmf}) in both the symmetric and nonsymmetric case. We present some simple algorithms in  Section \ref{s:algos}. Finally, in Section \ref{s:num-res}, we compare our methods with existing algorithms in numerical examples.

\section{Convex approximations}
\label{s:cvx-approx} In this section we derive a convex approximation to problem (\ref{eq:nmf}). We first discuss the case where the matrix $A$ is symmetric, then generalize our results to nonsymmetric matrices.

\subsection{Symmetric case}
\label{s:sym} In this section, given a data matrix $A\in\symm^n$, we focus on the following symmetric nonnegative matrix factorization problem:
\BEQ\label{eq:snnmf}
\BA{ll}
\mbox{minimize} & \mathbf{loss}(A,UU^T)\\
\mbox{subject to} & U \geq 0,\\
\EA\EEQ
in the variable $U\in\reals^{n \times k}$. In this paper, we consider two classical choices for the loss function given by:
\BC\begin{tabular}{rl}
Mean Squared Error (MSE): & $\mathbf{loss}(X,Y)=\|X-Y\|^2$\\
Kullback-Leibler (KL): & $\mathbf{loss}(X,Y)=\sum_{i,j=1}^n\left(X_{ij}\log(X_{ij}/Y_{ij})+Y_{ij}-X_{ij}\right)$.\\
\end{tabular}\EC
Both losses are convex functions of either $X$ or $Y$ but are not jointly in convex in $(X,Y)$, which means (\ref{eq:nmf}) is a nonconvex problem.

\subsubsection{Completely positive matrices} The solutions to this symmetric NMF problem, i.e. the symmetric matrices $A\in\symm^n$ which can be written in the form:
\BEQ \label{eq:cp}
A=UU^T, \quad U\geq 0
\EEQ
where $U\in\reals^{n \times k}$, are called \emph{completely positive} (see \cite{Berm03} for a complete discussion). This can also be written $A=\sum_{i=1}^k u_iu_i^T$ with $u_i\geq 0$ and the smallest $k$ for which this representation holds is called the \emph{cp-rank} of $A$. This provides us with another interpretation of problem (\ref{eq:snnmf}): it seeks the closest completely positive matrix to a given a matrix $A$. It also illustrates the difficulty of this problem: the set of completely positive matrices forms a cone whose dual is the cone of copositive matrices and testing for the copositivity of a matrix is a well-known NP-Hard problem. Also, \cite{Bure06} shows that any continuous or binary nonconvex program can be written as a linear program over the cone of completely positive matrices. The fundamental result we use in this section is given by the following theorem.
\begin{theorem}[Th. 2.30 in \cite{Berm03}] \label{th:cp-sdp}
If $X\in\symm^n$ is positive semidefinite, then $\exp_H(X)$, the Hadamard (or componentwise) exponential of $X$, is completely positive.
\end{theorem}
This means that after a natural change of variables $A_{ij}=\exp(X_{ij})$, we get a sufficient, \emph{convex} condition on $X$ for representation (\ref{eq:cp}) to hold. This also shows that kernel matrices obtained by negative exponentiation of negative semidefinite distance matrices are completely positive, hence can be interpreted as linear kernels over nonnegative feature vectors.

\subsubsection{Convex restriction} The above result allows us to form a convex restriction of the symmetric NMF problem in (\ref{eq:snnmf}) as:
\BEQ\label{eq:psd-snnmf}
\BA{ll}
\mbox{minimize} & \mathbf{loss}(A,\exp_H(X))\\
\mbox{subject to} & X \succeq 0,\\
\EA\EEQ
in the variable $X\in\symm^n$. When the loss function is given by the KL divergence, this problem becomes:
\BEQ\label{eq:psd-snnmf-kl}
\BA{ll}
\mbox{minimize} & \sum_{i,j=1}^n A_{ij}(\log(A_{ij})-X_{ij})+\exp(X_{ij})-A_{ij}\\
\mbox{subject to} & X \succeq 0,\\
\EA\EEQ
which is a \emph{convex} optimization problem in the variable $X\in\symm^n$. When the loss is given by the MSE, the objective  function is not convex in $X$ unless we impose the additional constraint that $A_{ij}/2\leq\exp(X_{ij})$. Problem (\ref{eq:psd-snnmf}) then becomes:
\BEQ\label{eq:psd-snnmf-mse}
\BA{ll}
\mbox{minimize} & \sum_{i,j=1}^n (\exp(X_{ij})-A_{ij})^2\\
\mbox{subject to} & A_{ij}/2\leq\exp(X_{ij}),\quad i,j=1,\ldots,n\\
& X \succeq 0,\\
\EA\EEQ
which is also a convex problem in the variable $X\in\symm^n$.

\subsubsection{Factorizing ${\exp_H(X)}$ when $X$ is positive semidefinite}
We know from Theorem \ref{th:cp-sdp} that ${\exp_H(X)}$ is completely positive when $X$ is positive semidefinite so there is a matrix $L$ such that 
\[{\exp_H(X)}=UU^T,\]
where $U\in\reals^{n \times k}$. First, Carath\'eodory's theorem allows us to bound the size of $U$ (i.e. the cp-rank of $\exp_H(X)$), and Theorem 3.5 in \cite{Berm03} shows that we can get:
\BEQ\label{eq:cp-rank-bound}
k\leq \frac{r(r+1)}{2}-1,
\EEQ
where $r=\Rank(\exp_H(X))$. Also, the Hadamard (or componentwise) product of two completely positive matrices is completely positive: suppose $A=\sum_{i=1}^ka_ia_i^T$ and $B=\sum_{i=1}^l b_ib_i^T$ with $a_i,b_i \geq 0$, then:
\BEQ\label{eq:had-prod}
A\circ B=\sum_{i=1}^k \sum_{j=1}^l (a_i \circ b_j)(a_i \circ b_j)^T,
\EEQ
hence the Hadamard product $A \circ B$ is completely positive as the sum of nonnegative rank one matrices. Now, because $X$ is positive semidefinite, we can write: 
\[
\exp_H(X)=\prod_{i=1}^n \exp_H\left(\lambda_i x_ix_i^T\right),
\]
where the matrix product is understood componentwise and $\lambda_i \in \reals^n_+$, which means that $\exp_H(X)$ can be written as the Hadamard product of matrices of the type $\exp_H(vv^T)$. As in \cite{Berm03} Theorem 2.30, we let $M=\max_{i=1,\ldots,n}|v_i|$, then
\[
\exp_H(vv^T)_{ij}=\exp(v_iv_j)=\exp(-M^2+(M+x_i)(M+x_j)-M(x_i+x_j)), \]
so
\[
\exp_H(vv^T)=\exp(-M^2)\exp_H(yy^T)\circ zz^T,
\]
where $y=M\ones +v$ and $z=\exp_H(-Mv)$ are both nonnegative vectors. Because $y$ is nonnegative and
\[
\exp_H(yy^T)=\sum_{i=1}^\infty \frac{(yy^T)^{\circ i}}{i!},
\]
with $(yy^T)^{\circ k}_{ij}=(y_iy_j)^k$, the matrix $\exp_H(yy^T)$ is completely positive. Hence, as the Hadamard product of two completely positive matrices, $\exp_H(vv^T)$ is completely positive. To summarize, when $X$ is positive semidefinite, we factorize the matrix $\exp_H(X)$ as follows:\\
\line(1,0){395}
\vskip 0ex
\textbf{Factorizing ${\exp_H(X)}$}
\begin{enumerate} \itemsep 0ex
\item Compute the eigenvalue decomposition: $X=\sum_{i=1}^n \lambda_i x_ix_i^T$.
\item Decompose each factor, $\exp_H(v_iv_i^T)=\exp(-M^2)\exp_H(y_iy_i^T)\circ z_iz_i^T$ where $v_i=\sqrt{\lambda_i}x_i$ and $y_i,z_i$ are nonnegative vectors.
\item Approximate $\exp_H(y_iy_i^T)$ as $\sum_{i=1}^k {(yy^T)^{\circ k}}/{k!}$.
\item Collect all the terms above using the chain rule in (\ref{eq:had-prod}) to get ${\exp_H(X)}=UU^T$.
\end{enumerate}
\line(1,0){395}
\vskip 0ex
Without any further processing, the size of $U$ can quickly become very large. Because of the bound given in (\ref{eq:cp-rank-bound}), we know however that the number of columns of $U$ is bounded above by ${r(r+1)}/{2}-1$ where $r=\Rank(\exp_H(X))$ and we can use this result to simplify the decomposition, but this is numerically costly and typically unnecessary. In practice, the eigenvalues of $\exp_H(v_iv_i^T)$ are decreasing exponentially fast and we can use the fact that when $X\in\symm^n$ is a positive semidefinite matrix with nonnegative coefficients and $\Rank(X)=2$, then $X$ is completely positive. We then replace $\exp_H(v_iv_i^T)$ by a rank two approximation which, if it is nonnegative, means that the size $k$ of the matrix $U\in\reals^{n \times k}$ is given by $k=2^r$ where $r$ is the number of significant factors in $X$, which is typically small. The precision of that approximation is further studied in Section \ref{s:num-res}.

\subsection{Sparse decomposition}
\label{ss:sparse-decomposition} As suggested in \cite{Hoye04} (see also \cite{Heil06}), e.g. when the matrix $A$ itself is sparse, we look for a \emph{sparse} decomposition $A=UU^T$, i.e. a decomposition where the matrices $U$ are sparse. In that case, the change of variable $A=\exp(X_{ij})$ is not appropriate. We can however look directly for a low rank decomposition by exploiting the property detailed above that when $X\in\symm^n$ is a positive semidefinite matrix with nonnegative coefficients and $\Rank(X)=2$, then $X$ is completely positive. We then solve:
\BEQ\label{eq:low-rank}
\BA{ll}
\mbox{minimize} & \|A-X\|^2+\gamma |X|+\nu\Tr(X)\\
\mbox{subject to} & X \geq 0,~ X \succeq 0,\\
\EA\EEQ
in the variable $X\in\symm^n$, where $|X|=\sum_{i,j=1}^n |X_{ij}|$, $\gamma$ is a parameter controlling the sparsity of $X$ and $\nu$ is a penalty on its rank (see \cite{Boyd00} for details).

\subsection{Recursive decomposition} Because the condition in Theorem \ref{th:cp-sdp} is only sufficient, problems (\ref{eq:psd-snnmf-kl}) and (\ref{eq:psd-snnmf-mse}) only cover a subset of all the possible nonnegative factorizations of the data matrix $A$. To overcome this limitation, we can solve problem (\ref{eq:snnmf}) recursively, setting $A_0=A$ and 
\[
A_{k+1}=A_{k}-U_kU_k^T,
\] 
where $U_k$ is the solution to the factorization problem given $A_k$. To ensure, that the intermediate matrices $A_k$ remain nonnegative, we can simply add the (convex) constraint that $U_kU_k^T\leq A$ to problems (\ref{eq:psd-snnmf-kl}) and (\ref{eq:psd-snnmf-mse}).

\subsection{Nonsymmetric case} We now extend the results of the previous section to the symmetric case. Given a data matrix $A\in\reals^{m \times n}$, we write the nonnegative matrix factorization problem as:
\BEQ\label{eq:nsnmf}
\BA{ll}
\mbox{minimize} & \mathbf{loss}(A,UV^T)\\
\mbox{subject to} & U,V \geq 0,\\
\EA\EEQ
in the variables $U\in\reals^{m \times k}$ and $V\in\reals^{n \times k}$, where $\mathbf{loss}(X,Y)$ is one of the  loss functions given above. Because the matrix $A$ has a nonnegative factorization if and only if there are matrices $B,C\in\symm^n$ such that the symmetric block-matrix:
\[
\left(\BA{cc}
B & A\\
A^T & C\\
\EA\right)
\]
is completely positive. Of course any nonnegative matrix can be factorized as the product of two nonnegative matrices because $A=A\mathbf{I}$ is always a solution. Our objective here is to make this representation as parsimonious as possible and find a solution with minimum cp-rank. We can't minimize the cp-rank of the decomposition directly without making the problem nonconvex, however the bound in (\ref{eq:cp-rank-bound}) shows that we can use the rank of a matrix as a proxy for its cp-rank. We know from \cite{Boyd00} that when $X\in\reals^{m \times n}$, $\|X\|_*$, the trace norm of $X$ is the largest convex lower bound on $\Rank(X)$. We also know that $\|X\|_*\leq t$ if and only if there are symmetric matrices $Y,Z$ such that:
\[
\left(\BA{cc}
Y & X\\
X^T & Z\\
\EA\right)\succeq 0
\quad \mbox{and} \quad
\Tr(Y)+\Tr(Z)\leq 2t.
\]
The problem of finding a low cp-rank nonnegative factorization of a matrix $A\in\reals^{m \times n}$ can then be written:
\BEQ\label{eq:cvx-nsnmf}
\BA{ll}
\mbox{minimize} & \mathbf{loss}(A,X) + \gamma (\Tr(Y)+\Tr(Z))\\
\mbox{subject to} & \left(\BA{cc}
Y & X\\
X^T & Z\\
\EA\right) \mbox{completely positive,}\\
\EA\EEQ
in the variables $X\in\reals^{m \times n}$, $Y\in\symm^m$ and $Z\in\symm^n$, where $\gamma>0$ controls the rank of the solution. This is now a symmetric NMF problem and can be solved using the results detailed in the previous sections.

\section{Algorithms}
\label{s:algos}
The results in (\ref{eq:psd-snnmf-kl}) and (\ref{eq:psd-snnmf-mse}) show that the symmetric NMF problem can be approximated by convex problems for which efficient, globally convergent algorithms are available. Here, because our focus is on solving large-scale problems with a relatively low precision, we use simple first-order methods to solve the optimization problems detailed in Section \ref{s:cvx-approx}. 

\subsection{Projected gradient method}
The simplest of these algorithms is the projected gradient method. Suppose that we need to minimize a convex function $f(x)$ over a convex set $C$. The projected gradient algorithms works as follows:\\
\line(1,0){395}
\vskip 0ex
\textbf{Projected Gradient Method}
\begin{enumerate} \itemsep 0ex
\item Start from a point $x_0\in\reals^n$.
\item Compute $x_{k+1}=p_C(x_k+\nabla f(x_k))$, where $p_C(x)$ is the projection of $x$ on the set $C$. 
\item Repeat step 2 until precision target is reached.
\end{enumerate}
\line(1,0){395}

Applying this method to problem (\ref{eq:psd-snnmf-kl}) for example, to solve:
\[\BA{ll}
\mbox{minimize} & \sum_{i,j=1}^n A_{ij}(\log(A_{ij})-X_{ij})+\exp(X_{ij})-A_{ij}\\
\mbox{subject to} & X \succeq 0,\\
\EA\]
we get $p_C$ explicitly as $p_C(X)=X_+$, where the matrix $X_+$ is obtained by zeroing out the negative eigenvalues of the matrix $X$.

\subsection{Complexity}
Provided some smoothness assumptions on the objective function, the number of iterations of the generic projected gradient method grows as $O(1/\epsilon^2)$ where $\epsilon$ is the target precision.

\subsection{Convergence and duality gap}
The dual of problem (\ref{eq:psd-snnmf-kl}) is given by:
\BEQ\label{eq:dual}
\BA{ll}
\mbox{maxmimize} & \sum_{i,j=1}^n A_{ij}+Y_{ij} -(A_{ij}+Y_{ij})\log(A_{ij}+Y_{ij})\\
\mbox{subject to} & Y \succeq 0,\\
\EA\EEQ
in the variable $Y\in \symm^n$. The optimality conditions impose:
\[
Y=\exp_H(X)-A\succeq 0,
\]
where $X,Y$ are primal and dual solutions to (\ref{eq:psd-snnmf-kl}). This means that if $X$ is the current primal point, then $(\exp_H(X)-A)_+$ is a dual feasible point which can be used to compute a duality gap, measure the optimality of $X$ and track convergence.

\section{Numerical Results}
\label{s:num-res}

In this section, we test the performance of our algorithm for solving the symmetric NMF problem on graph partitioning problems.

\subsection{Graph partitioning using NMF}
Let $A\in\{0,1\}^{n\times n}$ be the adjacency matrix of a given graph $G$, with
\[
A_{ij}=\left\{\BA{l}
1~\mbox{if $(i,j)$ is an edge of $G$}\\
0~\mbox{otherwise.}\EA\right.
\]
Suppose that we want to partition this graph into $k$ clusters $C_k$ while minimizing the number of graph edges between clusters and maximizing the number of graph edges inside clusters. For a given partition $C$ of a graph with adjacency matrix $A\in\{0,1\}^{n\times n}$, the performance measure we use here is given by:
\BEQ\label{eq:perf}
\mathbf{perf}(C)=1-\frac{\#\left\{(i,j)~|~A_{ij}\neq\sum_{l=1}^k X_{il}X_{jl}\right\}}{n^2}
\EEQ
where $X\in\{0,1\}^{n\times k}$ is an \emph{indicator matrix} such that:
\[
X_{ik}=\left\{\BA{l}
1~\mbox{if node $i$ is in cluster $C_k$}\\
0~\mbox{otherwise,}\EA\right.
\]
which satisfies $X\ones=\ones$. The graph partitioning problem can then be formulated as:
\BEQ\label{eq:graph-part}
\BA{ll}
\mbox{minimize} & \|A-XX^T\|^2\\
\mbox{subject to} & X\ones=\ones,
\EA\EEQ
in the variables $X\in\{0,1\}^{n\times k}$. This is a (hard) binary optimization problem, which we relax into a symmetric nonnegative matrix factorization problem as in \cite{Ding05}, to get:
\BEQ\label{eq:graph-part-relax}
\BA{ll}
\mbox{minimize} & \|A-XX^T\|^2\\
\mbox{subject to} & X\geq 0,
\EA\EEQ
in the variable $X\in\reals^{n\times k}$, which can be solved using the algorithms detailed in Section \ref{s:algos}. We can then turn the solution $X$ of this relaxation into an indicator matrix by setting the maximum coefficient of each row to one and all the others to zero.

\subsection{Partitioning performance}
\begin{figure}[htbp]
\begin{center}
\psfrag{orig}[t][b]{Original adjacency matrix}
\psfrag{new}[t][b]{Clustered permutation}
\includegraphics[width=1.0\textwidth]{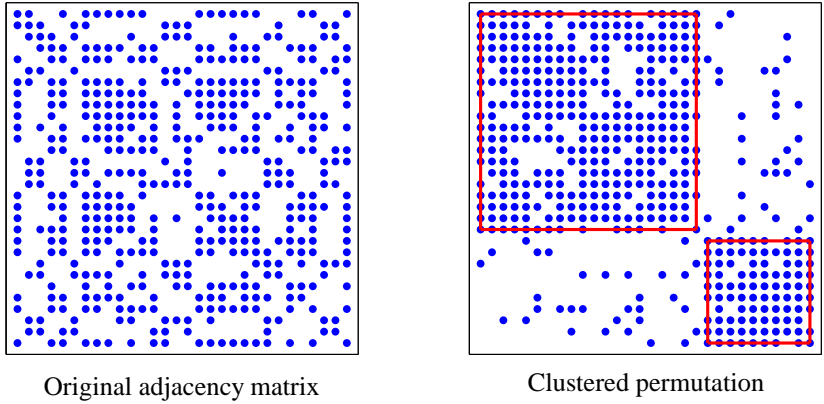}
\caption{Graph partitioning example: original (randomly generated) adjacency matrix on the left, clustered permutation on the right obtained by solving problem (\ref{eq:graph-part-relax}) using the algorithms in Section \ref{s:algos}.\label{fig:graph-part-ex}}
\end{center}
\end{figure}

To test the performance of our algorithms for symmetric NMF on graph partitioning problems, we generate random graphs whose adjacency matrices have given block sparsity patterns. We generate three uniform random matrices $A,B\in\symm^n$ and $C\in\reals^{n\times n}$ and use two parameters $\alpha,\beta\in[0,1]$ to control the sparsity and form a sample graph adjacency matrix as:
\[\left(\BA{cc}
1_{\{A_{ij}\geq \alpha\}} & 1_{\{C_{ij}\geq \beta\}}\\
1_{\{C_{ij}\geq \beta\}}^T & 1_{\{B_{ij}\geq \alpha\}}\\
\EA\right)\]
then randomly permute the matrix. An example is detailed in Figure \ref{fig:graph-part-ex}. 
We then compare the performance of our code with that of spectral clustering (see \cite{Ding01} for example) and show the results in Figure \ref{fig:stats} below. We observe that both methods perform similarly well on clearly clustered data but that the NMF solution dominates spectral clustering as the graphs become closer to bipartite.

\subsection{Convergence speed}
\begin{figure} [ht]
\psfrag{iter}[b][t]{\small{Number of iterations}}
\psfrag{n}[t][b]{\small{$n$}}
\includegraphics[width=.50 \textwidth]{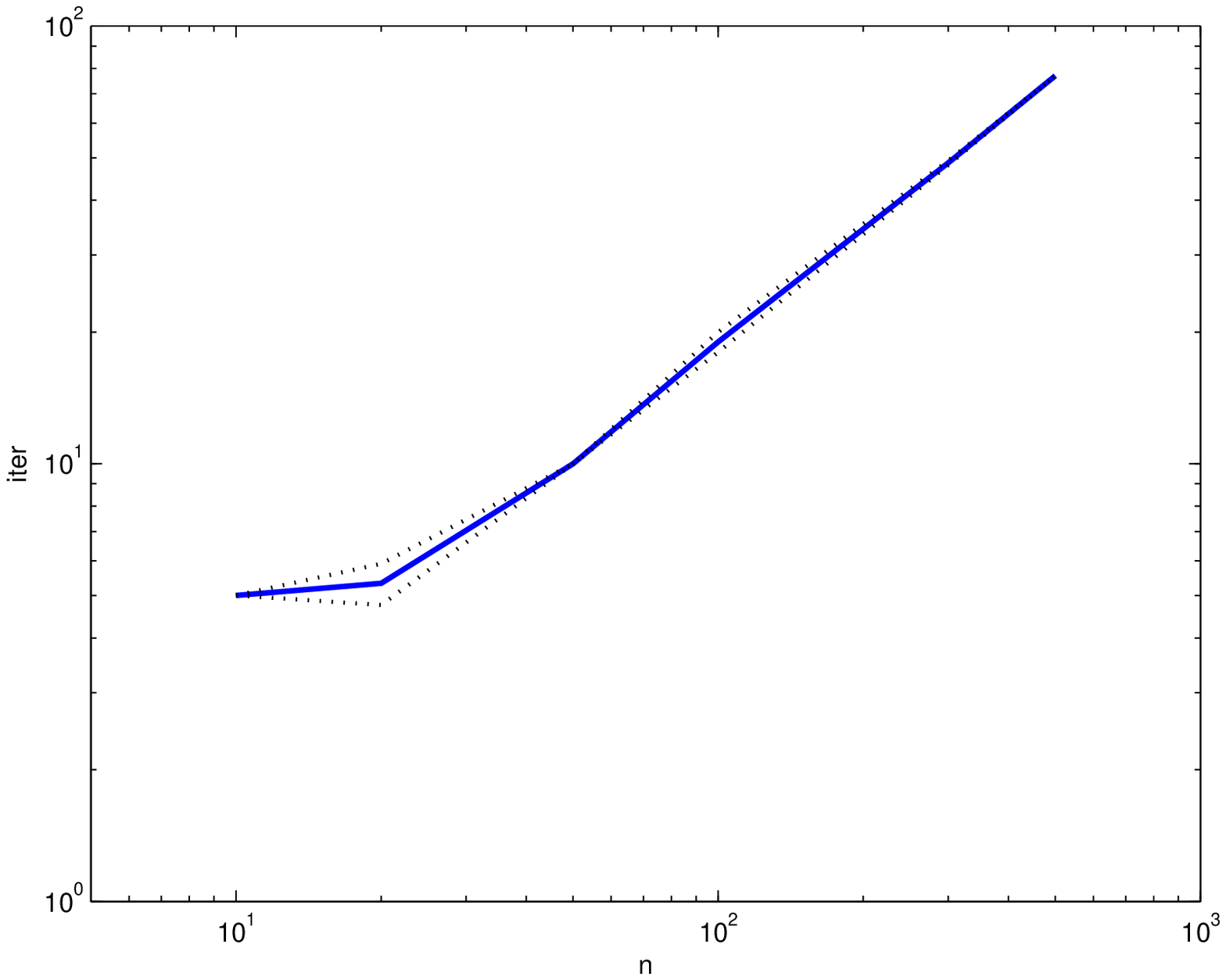} 
\psfrag{perf}[b][t]{\small{Performance}}
\psfrag{ab}[t][b]{\small{$\beta-\alpha$}}
\includegraphics[width=.50 \textwidth]{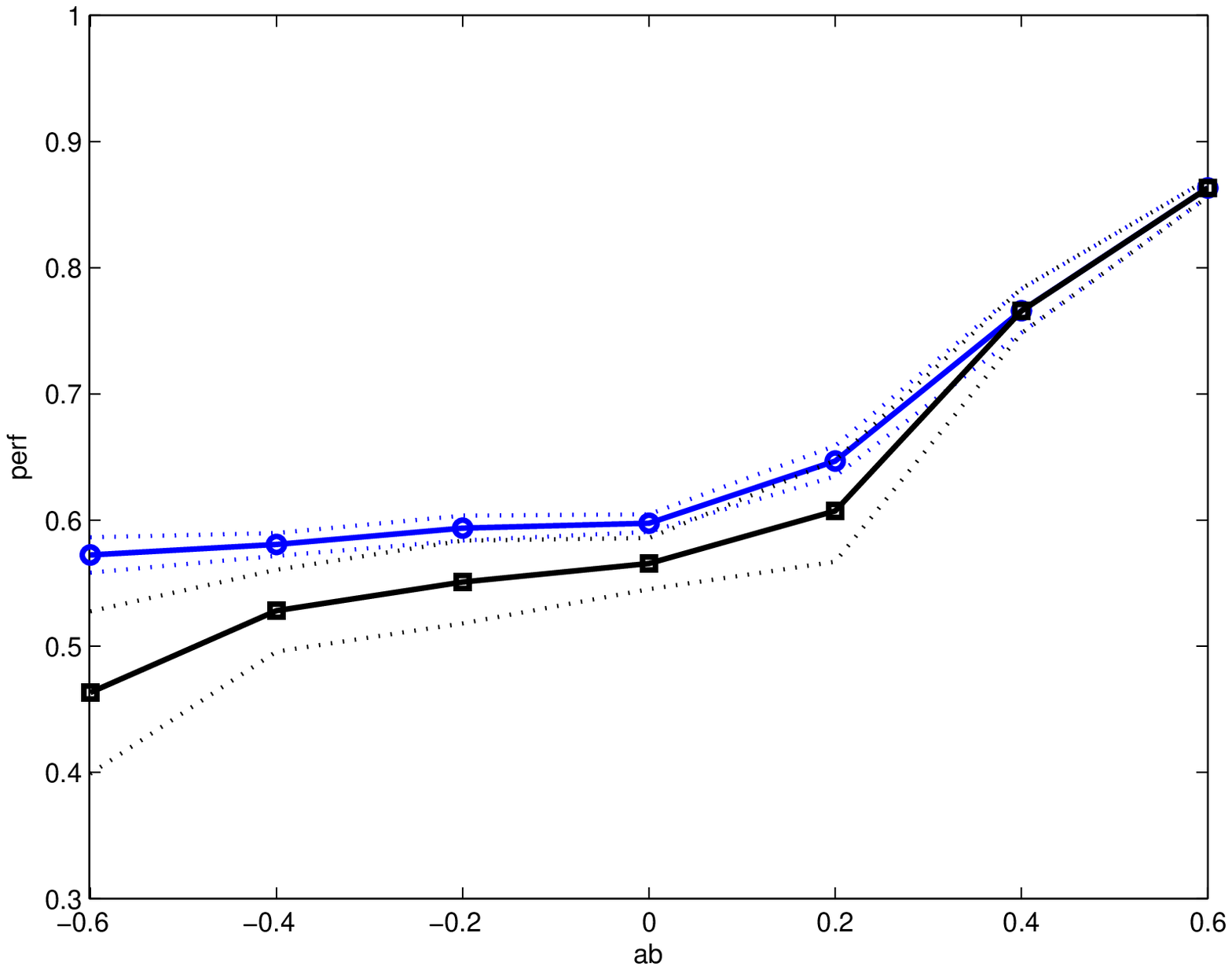}
\caption{\emph{Left}: Average number of iterations versus matrix size for the projected gradient algorithm. \emph{Right}: Average performance versus graph connectivity for spectral clustering (squares) and the solution to the NMF problem (\ref{eq:graph-part-relax}) using the algorithms in Section \ref{s:algos} (circles, dashed lines at plus and minus one standard deviation). \label{fig:stats}}
\end{figure}

In Figure \ref{fig:stats} we plot number of iterations of the projected gradient algorithm versus matrix size $n$ in randomly generated problems for various problem sizes.

\subsection{Extensions}
While the results in the symmetric case seem to perform very well in numerical experiments, this is not the case for nonsymmetric problems where alternating projection methods using highly optimized interior point solvers are still at least an order of magnitude faster than our method. At this point, efficiently exploiting the NMF representation detailed here in the nonsymmetric case remains an open numerical problem.

\bibliographystyle{unsrt}
\bibliography{MainPerso}

\begin{thebibliography}{10}

\bibitem{Paat94}
P.~Paatero and U.~Tapper.
\newblock Positive matrix factorization: A non-negative factor model with
  optimal utilization of error estimates of data values.
\newblock {\em Environmetrics}, 5(111-126), 1994.

\bibitem{Lee99}
D.D. Lee and H.S. Seung.
\newblock Learning the parts of objects by non-negative matrix factorization.
\newblock {\em Nature}, 401(6755):788--791, 1999.

\bibitem{Brun04}
J.P. Brunet, P.~Tamayo, T.R. Golub, and J.P. Mesirov.
\newblock Metagenes and molecular pattern discovery using matrix factorization.
\newblock {\em Proceedings of the National Academy of Sciences},
  101(12):4164--4169, 2004.

\bibitem{Sha05}
F.~Sha and L.K. Saul.
\newblock Real-time pitch determination of one or more voices by nonnegative
  matrix factorization.
\newblock {\em Advances in Neural Information Processing Systems}, 17, 2005.

\bibitem{Xu03}
Wei Xu, Xin Liu, and Yihong Gong.
\newblock Document clustering based on non-negative matrix factorization.
\newblock In {\em SIGIR '03: Proceedings of the 26th annual international ACM
  SIGIR conference on Research and development in informaion retrieval}, pages
  267--273, 2003.

\bibitem{Lang05a}
T.~Lange and J.M. Buhmann.
\newblock Fusion of similarity data in clustering.
\newblock {\em Advances in Neural Information Processing Systems 18}, 2005.

\bibitem{Dono03}
D.~Donoho and V.~Stodden.
\newblock When does non-negative matrix factorization give a correct
  decomposition into parts?
\newblock {\em Advances in Neural Information Processing Systems 15}, 2003.

\bibitem{Berr06}
M.W. Berry, M.~Browne, A.N. Langville, V.P. Pauca, and R.J. Plemmons.
\newblock Algorithms and applications for approximate nonnegative matrix
  factorization.
\newblock {\em To appear in Computational Statistics and Data Analysis}, 2006.

\bibitem{Zdun06}
R.~Zdunek and A.~Cichocki.
\newblock Non-negative matrix factorization with quasi-newton optimization.
\newblock {\em Eighth International Conference on Artificial Intelligence and
  Soft Computing, ICAISC}, pages 870--879, 2006.

\bibitem{Hoye04}
P.O. Hoyer.
\newblock Non-negative matrix factorization with sparseness constraints.
\newblock {\em The Journal of Machine Learning Research}, 5:1457--1469, 2004.

\bibitem{Lee01}
D.D. Lee and H.S. Seung.
\newblock Algorithms for non-negative matrix factorization.
\newblock {\em Advances in Neural Information Processing Systems}, 13:556--562,
  2001.

\bibitem{Heil06}
M.~Heiler and C.~Schnorr.
\newblock Learning sparse representations by non-negative matrix factorization
  and sequential cone programming.
\newblock {\em Journal of Machine Learning Research}, 7:1385--1407, 2006.

\bibitem{Berm03}
A.~Berman and N.~Shaked-Monderer.
\newblock {\em Completely Positive Matrices}.
\newblock World Scientific, 2003.

\bibitem{Bure06}
S.~Burer.
\newblock On the copositive representation of binary and continuous nonconvex
  quadratic programs.
\newblock {\em Working Paper}, 2006.

\bibitem{Boyd00}
M.~Fazel, H.~Hindi, and S.~Boyd.
\newblock A rank minimization heuristic with application to minimum order
  system approximation.
\newblock {\em Proceedings American Control Conference}, 6:4734--4739, 2001.

\bibitem{Ding05}
C.~Ding, X.~He, and H.~Simon.
\newblock On the equivalence of nonnegative matrix factorization and spectral
  clustering.
\newblock {\em Proceedings of the SIAM Data Mining Conference}, 2005.

\bibitem{Ding01}
C.~Ding, X.~He, H.~Zha, M.~Gu, and HD~Simon.
\newblock A min-max cut algorithm for graph partitioning and data clustering.
\newblock {\em Proceedings IEEE International Conference on Data Mining}, pages
  107--114, 2001.

\end{thebibliography}
\end{document}